\documentclass[12pt, twoside, leqno]{article}
\usepackage{amsmath,amsthm}
\usepackage{amssymb,latexsym}
\usepackage{enumerate}

\newtheorem{thm}{Theorem}[section]
\newtheorem{prop}[thm]{Proposition}
\newtheorem{cor}[thm]{Corollary}
\newtheorem{lem}[thm]{Lemma}

\newtheorem{conj}[thm]{Conjecture}

\theoremstyle{definition}

\newtheorem{rem}[thm]{Remark}

\numberwithin{equation}{section}

\usepackage{mathrsfs}
\usepackage{CJK}
\usepackage{amsmath, amscd}
\usepackage{fancyhdr}
\usepackage{enumerate}

\begin{document}

\baselineskip=17pt

\title{The LFED Conjecture for some $\mathcal{E}$-derivations}
\author{ Lintong Lv \\
MOE-LCSM,\\ School of Mathematics and Statistics,\\
 Hunan Normal University, Changsha 410081, China \\
 \emph{E-mail:} lvlintong97@163.com \\
\and
Dan Yan \footnote{ The author is supported by the NSF of China (Grant No. 11871241; 11601146), the China Scholarship Council and the Construct Program of the Key Discipline in Hunan Province.}\\
MOE-LCSM,\\ School of Mathematics and Statistics,\\
 Hunan Normal University, Changsha 410081, China \\
\emph{E-mail:} yan-dan-hi@163.com
}
\date{}

\maketitle

\renewcommand{\thefootnote}{}

\renewcommand{\thefootnote}{\arabic{footnote}}
\setcounter{footnote}{0}

\begin{abstract}
 Let $K$ be an algebraically closed field of characteristic zero, $\delta$ a nonzero $\mathcal{E}$-derivation of $K[x]$. We first prove that $\operatorname{Im}\delta$ is a Mathieu-Zhao space of $K[x]$ in some cases. Then we prove that LFED Conjecture is true for all $\delta=I-\phi$, where $\phi$ is an affine polynomial homomorphism of $K[x_1,x_2]$. Finally, we prove that LFED Conjecture is true for some $\delta$ of $K[x_1,x_2,x_3]$.
\end{abstract}
{\bf Keywords.} Mathieu-Zhao spaces, $\mathcal{E}$-derivations, Locally finite.\\
{\bf MSC(2010).} 13N15; 14R10; 13F20. \vskip 2.5mm

\section{Introduction}

Throughout this paper, we will write $K$ for an algebraically closed field of characteristic zero without specific note and $K[x]=K[x_1,x_2,\ldots,x_n]$ for the
polynomial algebra over $K$ with $n$ indeterminates. $\partial_i$ denotes the derivations $\frac{\partial}{\partial x_i}$ for $1\leq i\leq n$.

A $K$-linear endomorphism $\eta$ of $K[x]$ is said to be locally nilpotent if for each $a\in K[x]$ there exists $m\geq 1$ such that $\eta^m(a)=0$, and locally finite if for each $a\in K[x]$ the $K$-subspace spanned by $\eta^i(a)$ ($i\geq 0$) over $K$ is finitely generated.

A derivation $D$ of $K[x]$ we mean a $K$-linear map $D: K[x]\rightarrow K[x]$ that satisfies $D(ab)=D(a)b+aD(b)$ for all $a, b\in K[x]$. An $\mathcal{E}$-derivation $\delta$ of $K[x]$ we mean a $K$-linear map $\delta: K[x]\rightarrow K[x]$ such that for all $a, b\in K[x]$ the following equation holds:
$$\delta(ab)=\delta(a)b+a\delta(b)-\delta(a)\delta(b).$$

It is easy to verify that $\delta$ is an $\mathcal{E}$-derivation of $K[x]$, if and only if $\delta=I-\phi$ for some $K$-algebra endomorphism $\phi$ of $K[x]$.

The Mathieu-Zhao space was introduced by Zhao in \cite{1} and \cite{2}, which is a natural generalization of ideals. We give the definition here for the polynomial rings. A $K$-subspace $M$ of $K[x]$ is said to be a Mathieu-Zhao space if for any $a, b\in K[x]$ with $a^m \in M$ for all $m\geq 1$, we have $ba^m\in M$ when $m>>0$. The radical of a Mathieu-Zhao space was first introduced in \cite{2}, denoted by $\mathfrak{r}(M)$, and
$$\mathfrak{r}(M)=\{a\in K[x]| a^m\in M~\operatorname{for}~\operatorname{all}~m>>0\}.$$

There is an equivalent definition about Mathieu-Zhao space which proved in Proposition 2.1 of \cite{2}. We only give the equivalent definition here for the polynomial rings. A $K$-subspace $M$ of $K[x]$ is said to be a Mathieu-Zhao space if for any $a, b\in K[x]$ with $a\in \mathfrak{r}(M)$, we have $ba^m\in M$ when $m>>0$.

In \cite{3}, Wenhua Zhao posed the following two conjectures:
\begin{conj}(LFED)\label{conj1.1}
Let $K$ be a field of characteristic zero and $\mathcal{A}$ a $K$-algebra. Then for every locally finite derivation or $\mathcal{E}$-derivation $\delta$ of $\mathcal{A}$, the image $\operatorname{Im}\delta:=\delta(\mathcal{A})$ of $\delta$ is a Mathieu-Zhao space of $\mathcal{A}$.
\end{conj}

\begin{conj}(LNED)\label{conj1.2}
Let $K$ be a field of characteristic zero and $\mathcal{A}$ a $K$-algebra and $\delta$ a locally nilpotent derivation or $\mathcal{E}$-derivation of $\mathcal{A}$, Then for every $\vartheta$-ideal $I$ of $\mathcal{A}$, the image $\delta(I)$ of $I$ under $\delta$ is a $\vartheta$-MZ space of $\mathcal{A}$.
\end{conj}

There are many positive answers to the above two conjectures.  In \cite{4}, Wenhua Zhao proved that Conjecture \ref{conj1.1} is true for polynomial algebras in one variable and Conjecture \ref{conj1.2} is true for polynomial algebras in one variable for derivations and most $\mathcal{E}$-derivations. Arno van den Essen, David Wright, Wenhua Zhao showed that Conjecture \ref{conj1.1} is true for derivations for polynomial algebras in two variables in \cite{5}. In \cite{6}, Wenhua Zhao proved that Conjecture \ref{conj1.1} is true for Laurent polynomial algebras in one or two variables and Conjecture \ref{conj1.2} is true for all Laurent polynomial algebras. Wenhua Zhao proved the above two conjectures for algebraic algebras in \cite{7}. In \cite{8}, Dayan Liu, Xiaosong Sun showed that Conjecture \ref{conj1.1} is true for linear locally nilpotent derivations in dimension three. Arno van den Essen, Wenhua Zhao showed that Conjecture \ref{conj1.1} is true for locally integral domains and $K[[x]][x^{-1}]$ in \cite{9}.

In our paper, we prove that Conjecture \ref{conj1.1} is true for some derivations and $\mathcal{E}$-derivations of $K[x]$ in section 2. In section 3, we show that Conjecture \ref{conj1.1} is true for all $\mathcal{E}$-derivations $\delta=I-\phi$, where $\phi$ is an affine polynomial homomorphism of $K[x_1,x_2]$. Then we prove that Conjecture \ref{conj1.1} is true for most $\mathcal{E}$-derivations $\delta=I-\phi$, where $\phi$ is a linear polynomial homomorphism of $K[x_1,x_2,x_3]$ and give a conjecture for other $\mathcal{E}$-derivations $\delta=I-\phi$, where $\phi$ is a linear polynomial homomorphism of $K[x_1,x_2,x_3]$ in section 4.

\section{Conjecture \ref{conj1.1} for some derivations and $\mathcal{E}$-derivations}

\begin{thm} \label{thm2.1}
Let $\delta=I-\phi$ be an $\mathcal{E}$-derivation of $K[x]$ and $\phi=Ax$ is a linear polynomial homomorphism of $K[x]$ with $A\in \operatorname{M}_n(K)$. If $\lambda_{kk}^{i_k}\cdots \lambda_{nn}^{i_n}\neq 1$ for all $i_k,\ldots,i_n\in \mathbb{N}$, $i_k+\cdots+i_n\geq 1$, $1\leq k\leq n$, where $\lambda_{11},\ldots,\lambda_{nn}$ are the eigenvalues of $A$, then $\operatorname{Im}\delta$ is an ideal of $K[x]$. In particular, if $\lambda_{11}=\cdots=\lambda_{nn}:=\lambda$, then $\operatorname{Im}\delta$ is an ideal of $K[x]$ in the case that $\lambda$ is not a root of unity.
\end{thm}
\begin{proof}
Since $\phi=Ax$, there exists $T\in \operatorname{GL}_n(K)$ such that

$$T^{-1}AT=\left( \begin{matrix}
\lambda_{11} & \lambda_{12} &\cdots & \lambda_{1n} \\
0 & \lambda_{22} & \cdots & \lambda_{2n} \\
\vdots & \vdots &\ddots& \vdots\\
0 & 0 &\cdots & \lambda_{nn}
\end{matrix} \right).$$
Let $\sigma(x)=Tx$. Then we have $\sigma^{-1} \delta \sigma=I-\sigma^{-1} \phi \sigma$. It suffices to prove that $\operatorname{Im}(\sigma^{-1} \delta \sigma)$ is an ideal of $K[x]$. Let $\tilde{\delta}=\sigma^{-1} \delta \sigma=I-\tilde{\phi}$. Then $\tilde{\phi}(x_i)=\sum_{j=i}^n\lambda_{ij}x_j$ for $1\leq i\leq n$. Thus, we have
$$\tilde{\delta}(x_n^{i_n})=(1-\lambda_{nn}^{i_n})x_n^{i_n}.$$
Since $\lambda_{nn}^{i_n}\neq 1$, we have $x_n^{i_n}\in \operatorname{Im}\tilde{\delta}$ for all $i_n\in \mathbb{N}^*$.
 Suppose that $x_k^{l_k}x_{k+1}^{l_{k+1}}\cdots x_{n-1}^{l_{n-1}}x_n^{i_n}\allowbreak\in \operatorname{Im}\tilde{\delta}$ for $l_k\leq i_k-1$ or $l_k+l_{k+1}\leq i_k+i_{k+1}-1$ or $\cdots$ or $l_k+\cdots+l_{n-1}\leq i_k+\cdots+i_{n-1}-1$, $l_k+\cdots+l_{n-1}+i_n\geq 1$. Then we have
\begin{eqnarray}
\nonumber
\tilde{\delta}(x_k^{i_k}x_{k+1}^{i_{k+1}}\cdots x_n^{i_n})= x_k^{i_k}\cdots x_n^{i_n}-(\lambda_{kk}x_k+\cdots+\lambda_{kn}x_n)^{i_k}\cdots (\lambda_{nn}x_n)^{i_n}\\
\nonumber
                    =(1-\lambda_{kk}^{i_k}\cdots \lambda_{nn}^{i_n})x_k^{i_k}\cdots                     x_n^{i_n}+\operatorname{polynomial}~\operatorname{in}~\operatorname{Im}\tilde{\delta}.\nonumber
\end{eqnarray}
Since $\lambda_{kk}^{i_k}\cdots \lambda_{nn}^{i_n}\neq 1$, we have $x_k^{i_k}\cdots                     x_n^{i_n}\in \operatorname{Im}\tilde{\delta}$ for all $i_k,\ldots,i_n\in \mathbb{N}$, $i_k+\cdots+i_n\geq 1$, $1\leq k\leq n$. Since $1\notin\operatorname{Im}\tilde{\delta}$, we have that $\operatorname{Im}\tilde{\delta}$ is an ideal generated by $x_1,x_2,\ldots,x_n$. Then the conclusion follows.
\end{proof}

\begin{lem}\label{lem2.2}
Let $\delta=I-\phi$ be an $\mathcal{E}$-derivation of $K[x]$ and $\phi(x_i)=\lambda_ix_i+f_i(x_{i+1},\ldots,x_n)$ with $f_i(x_{i+1},\ldots,x_n)\in K[x_{i+1},\ldots,x_n]$ and $f_n\in K$ for $1\leq i\leq n-1$. If $\lambda_i\neq 1$, then there exists $\sigma\in \operatorname{Aut}(K[x])$ such that $\sigma^{-1}\delta\sigma=I-\tilde{\phi}$ and $\tilde{\phi}(x_i)=\lambda_ix_i+\tilde{f}_i(x_{i+1},\ldots,x_n)$, where $\tilde{f}_i(0)=0$ for all $1\leq i\leq n$.
\end{lem}
\begin{proof}
Let $\tilde{\delta}=\sigma^{-1}\delta\sigma$ and $\sigma(x_i)=x_i+c_i$ for all $1\leq i\leq n$, where $c_i=(\lambda_i-1)^{-1}f_i(-c_{i+1},\ldots,-c_n)$ and $c_n=(\lambda_n-1)^{-1}f_n$ for all $1\leq i\leq n-1$. Then $\tilde{\delta}=\sigma^{-1}\delta\sigma=I-\tilde{\phi}$ and
$$\tilde{\phi}(x_i)=\lambda_ix_i+(1-\lambda_i)c_i+f_i(x_{i+1}-c_{i+1},\ldots,x_n-c_n)$$
for all $1\leq i\leq n$. Let $\tilde{f}_i(x_{i+1},\ldots,x_n)=f_i(x_{i+1}-c_{i+1},\ldots,x_n-c_n)+(1-\lambda_i)c_i$ for $1\leq i\leq n$. Then the conclusion follows.
\end{proof}

\begin{prop}\label{prop2.3}
Let $\delta=I-\phi$ be an $\mathcal{E}$-derivation of $K[x]$ and $\phi(x_i)=\lambda_ix_i+f_i(x_{i+1},\ldots,x_n)$, where $f_i\in K[x_{i+1},\ldots,x_n]$ and $\lambda_i\in K$ for all $1\leq i\leq n$. If $\lambda_k^{i_k}\cdots \lambda_n^{i_n}\neq 1$ for all $i_k,\ldots,i_n\in \mathbb{N}$, $i_k+\cdots+i_n\geq 1$, $1\leq k\leq n$, then $\operatorname{Im}\delta$ is an ideal generated by $x_1,x_2,\ldots,x_n$.
\end{prop}
\begin{proof}
It follows from Lemma \ref{lem2.2} that we can assume that $f_i(0)=0$ for all $1\leq i\leq n$. Thus, we have $f_n=0$. Since $\phi(x_n)=\lambda_nx_n$, we have $\delta(x_n^{i_n})=(1-\lambda_n^{i_n})x_n^{i_n}$. Since $\lambda_n^{i_n}\neq 1$, we have $x_n^{i_n}\in \operatorname{Im}\delta$ for all $i_n\in \mathbb{N}^*$. Suppose that $x_k^{l_k}x_{k+1}^{l_{k+1}}\cdots x_{n-1}^{l_{n-1}}x_n^{i_n}\in \operatorname{Im}\delta$ for $l_k\leq i_k-1$ or $l_k+l_{k+1}\leq i_k+i_{k+1}-1$ or $\cdots$ or $l_k+\cdots+l_{n-1}\leq i_k+\cdots+i_{n-1}-1$, $l_k+\cdots+l_{n-1}+i_n\geq 1$. Then we have
\begin{eqnarray}
\nonumber
\delta(x_k^{i_k}x_{k+1}^{i_{k+1}}\cdots x_n^{i_n})= x_k^{i_k}\cdots x_n^{i_n}-(\lambda_kx_k+f_k(x_{k+1},\ldots,x_n))^{i_k}\cdots (\lambda_nx_n)^{i_n}\\
\nonumber
                    =(1-\lambda_k^{i_k}\cdots \lambda_n^{i_n})x_k^{i_k}\cdots                     x_n^{i_n}+\operatorname{polynomial}~\operatorname{in}~\operatorname{Im}\delta.\nonumber
\end{eqnarray}
Since $\lambda_k^{i_k}\cdots \lambda_n^{i_n}\neq 1$, we have $x_k^{i_k}\cdots                     x_n^{i_n}\in \operatorname{Im}\delta$ for all $i_k,\ldots,i_n\in \mathbb{N}$, $i_k+\cdots+i_n\geq 1$, $1\leq k\leq n$. Since $1\notin \operatorname{Im}\delta$, we have that $\operatorname{Im}\delta$ is an ideal generated by $x_1, x_2,\ldots,x_n$.
\end{proof}

\begin{prop}\label{prop2.4}
Let $\delta=I-\phi$ be an $\mathcal{E}$-derivation of $K[x]$ and $\phi(x_i)=\lambda_ix_i+\mu_i$, where $\lambda_i, \mu_i\in K$ for all $1\leq i\leq n$. Then $\operatorname{Im}\delta$ is a Mathieu-Zhao space of $K[x]$.
\end{prop}
\begin{proof}
If $\lambda_i\neq 1$ for some $i\in \{1,2,\ldots,n\}$, then we have $\sigma_i^{-1}\phi\sigma_i(x_i)=\lambda_ix_i$ and $\sigma_i^{-1}\phi\sigma_i(x_j)=\lambda_jx_j+\mu_j$, where $\sigma_i(x_i)=x_i+(\lambda_i-1)^{-1}\mu_i$, $\sigma_i(x_j)=x_j$ for $j\neq i$ for all $1\leq j\leq n$.

If $\lambda_i=1$, then $\delta(x_i)=-\mu_i$. If $\mu_i\neq 0$, then $1\in \operatorname{Im}\delta$. It's easy to check that $\delta$ is locally finite, it follows from Proposition 1.4 in \cite{10} $\operatorname{Im}\delta$ is a Mathieu-Zhao space of $K[x]$. If $\mu_i=0$, then $\phi(x_i)=\lambda_ix_i$. We assume that $\sigma_i=I$ in this case. Let $\sigma=\sigma_n\circ\cdots\circ\sigma_1\in \operatorname{Aut}(K[x])$. Then $\sigma^{-1}\delta\sigma=I-\tilde{\phi}$, where $\tilde{\phi}(x_i)=\lambda_ix_i$ for all $1\leq i\leq n$ or $\operatorname{Im}\delta$ is a Mathieu-Zhao space of $K[x]$. Let $\tilde{\delta}=\sigma^{-1}\delta\sigma$. It follows from Lemma 3.2 and Corollary 3.3 in \cite{11} that $\operatorname{Im}\tilde{\delta}$ is a Mathieu-Zhao space of $K[x]$. Thus, $\operatorname{Im}\delta$ is a Mathieu-Zhao space of $K[x]$.
\end{proof}

\begin{prop}\label{prop2.5}
Let $D=\sum_{i=1}^n(a_ix_i+b_i)\partial_i$ be a derivation of $K[x]$ with $a_i, b_i\in K$ for all $1\leq i\leq n$. Then $\operatorname{Im}D$ is a Mathieu-Zhao space of $K[x]$.
\end{prop}
\begin{proof}
If $a_i\neq 0$ for some $i\in \{1,2,\ldots,n\}$, then we have $$\ \sigma_i^{-1}D\sigma_i=a_ix_i\partial_i+\sum_{\mbox {\tiny $\begin{array}{c}
1\leq j\leq n\\
j\neq i\end{array}$}}(a_jx_j+b_j)\partial_j,$$
where $\sigma_i(x_i)=a_ix_i+b_i$, $\sigma_i(x_j)=x_j$ for $j\neq i$ for all $1\leq j\leq n$.

If $a_i=0$, then $D(x_i)=b_i$. If $b_i\neq 0$, then $1\in \operatorname{Im}D$. It follows from Example 9.3.2 in \cite{12} that $D$ is locally finite. Thus, it follows from Proposition 1.4 in \cite{10} that $\operatorname{Im}D$ is a Mathieu-Zhao space of $K[x]$. If $b_i=0$, then
$$D=\sum_{\mbox {\tiny $\begin{array}{c}
1\leq j\leq n\\
j\neq i\end{array}$}}(a_jx_j+b_j)\partial_j.$$
Hence we have that $\operatorname{Im}D$ is a Mathieu-Zhao space of $K[x]$ or there exists $\sigma\in \operatorname{Aut}(K[x])$ such that $\sigma^{-1}D\sigma=\sum_{j=1}^na_jx_j\partial_j$. It follows from Lemma 3.4 in \cite{5} that $\operatorname{Im}(\sigma^{-1}D\sigma)$ is a Mathieu-Zhao space of $K[x]$. Thus, $\operatorname{Im}D$ is a Mathieu-Zhao space of $K[x]$.
\end{proof}

\begin{prop}\label{prop2.6}
 Let $D=\sum_{i=1}^n(a_ix_i+b_i(x_1,\ldots,x_{i-1}))\partial_i$ be a derivation of $K[x]$ with $a_i\in K$, $b_i\in K[x_1,\ldots,x_{i-1}]$ for all $1\leq i\leq n$ and $S$ the set of nonzero integral solutions of the linear equation $\sum_{i=1}^na_iy_i=0$. If $S=\emptyset$, then $\operatorname{Im}D$ is an ideal of $K[x]$.
\end{prop}
\begin{proof}
Since $S=\emptyset$, we have $a_1a_2\cdots a_n\neq 0$. Thus, we have
$$\sigma_1^{-1}D\sigma_1=a_1x_1\partial_1+\sum_{i=2}^n(a_ix_i+b_i^{(1)}(x_1,\ldots,x_{i-1}))\partial_i,$$ where $\sigma_1(x_1)=a_1x_1+C_1$, $\sigma_1(x_i)=x_i$ and $b_i^{(1)}(x_1,\ldots,x_{i-1})=b_i(a_1^{-1}(x_1-C_1),x_2,\ldots,x_{i-1})$ for $2\leq i\leq n$ and $C_1=b_1$. Suppose that there exist polynomial automorphisms $\sigma_2,\ldots,\sigma_{k-1}$ such that
$$D_{k-1}:=\sigma_{k-1}^{-1}\cdots\sigma_1^{-1}D\sigma_1\cdots\sigma_{k-1}=\sum_{i=1}^{k-1}a_ix_i\partial_i+\sum_{j=k}^n(a_jx_j+b_j^{(k-1)}(x_1,\ldots,x_{j-1}))\partial_j.$$
We claim that there is a polynomial automorphism $\sigma_k$ such that
$$D_k=\sigma_k^{-1}D_{k-1}\sigma_k=\sum_{i=1}^ka_ix_i\partial_i+\sum_{j=k+1}^n(a_jx_j+b_j^{(k)}(x_1,\ldots,x_{j-1}))\partial_j.$$
Let $\sigma_k(x_i)=x_i$ for $i\neq k$ and for all $1\leq i\leq n$, $\sigma_k(x_k)=a_kx_k+C_k$ for $C_k\in K[x_1,\ldots,x_{k-1}]$. Then we have $\sigma_k^{-1}D_{k-1}\sigma_k(x_i)=a_ix_i$ for all $1\leq i\leq k-1$ and $\sigma_k^{-1}D_{k-1}\sigma_k(x_k)=a_kx_k-a_kC_k+a_kb_k^{(k-1)}(x_1,\ldots,x_{k-1})+\sum_{i=1}^{k-1}a_ix_i\frac{\partial C_k}{\partial x_i}$ and $\sigma_k^{-1}D_{k-1}\sigma_k(x_j)=a_jx_j+b_j^{(k-1)}(x_1,\ldots,x_{k-1},a_k^{-1}(x_k-C_k),x_{k+1},\ldots,x_{j-1})$
for all $k+1\leq j\leq n$.

Let $b_j^{(k)}:=b_j^{(k-1)}(x_1,\ldots,x_{k-1},a_k^{-1}(x_k-C_k),\ldots,x_{j-1})$ for all $k+1\leq j\leq n$. Then it suffices to prove that there exists $C_k\in K[x_1,\ldots,x_{k-1}]$ such that
\begin{eqnarray}\label{eq2.1}
a_kC_k-\sum_{i=1}^{k-1}a_ix_i\frac{\partial C_k}{\partial x_i}=a_kb_k^{(k-1)}(x_1,\ldots,x_{k-1}).
\end{eqnarray}
Let
$$b_k^{(k-1)}:=b_k^{(k-1)}(x_1,\ldots,x_{k-1})=\sum_{l_1,\ldots,l_{k-1}\geq 0}b_{k,l_1,\ldots,l_{k-1}}^{(k-1)}x_1^{l_1}\cdots x_{k-1}^{l_{k-1}}$$
and
$$C_k=\sum_{l_1,\ldots,l_{k-1}\geq 0}C_{k,l_1,\ldots,l_{k-1}}^{(k-1)}x_1^{l_1}\cdots x_{k-1}^{l_{k-1}}.$$
Then we have
$$a_kC_k-\sum_{i=1}^{k-1}a_ix_i\frac{\partial C_k}{\partial x_i}=\sum_{l_1,\ldots,l_{k-1}\geq 0}(a_k-\sum_{i=1}^{k-1}l_ia_i)C_{k,l_1,\ldots,l_{k-1}}x_1^{l_1}\cdots x_i^{l_i}\cdots x_{k-1}^{l_{k-1}}.$$
Thus, equation \eqref{eq2.1} is equivalent to that $(a_k-\sum_{i=1}^{k-1}l_ia_i)C_{k,l_1,\ldots,l_{k-1}}=b_{k,l_1,\ldots,l_{k-1}}^{(k-1)}a_k$.
Since $S=\emptyset$, we have $a_k-\sum_{i=1}^{k-1}l_ia_i\neq 0$ for all $l_1,\ldots,l_{k-1}\in \mathbb{N}$. Let $C_{k,l_1,\ldots,l_{k-1}}=(a_k-\sum_{i=1}^{k-1}l_ia_i)^{-1}a_kb_{k,l_1,\ldots,l_{k-1}}^{(k-1)}$ for all $l_1,\ldots,l_{k-1}\in \mathbb{N}$. Then we have $\sigma_k^{-1}D_{k-1}\sigma_k(x_k)\allowbreak=a_kx_k$. This finished the proof of our claim. Hence there exists $\sigma\in \operatorname{Aut}(K[x])$ such that $\sigma^{-1}D\sigma=\sum_{i=1}^na_ix_i\partial_i$.
Since $\sigma^{-1}D\sigma(x_1^{j_1}\cdots x_n^{j_n})=(\sum_{i=1}^nj_ia_i)x_1^{j_1}\cdots x_n^{j_n}$ for all $j_1,\ldots,j_n\in \mathbb{N}$, $j_1+\cdots+j_n\geq 1$ and $\sum_{i=1}^nj_ia_i\neq 0$, we have that $x_1^{j_1}\cdots x_n^{j_n}\in \operatorname{Im}\sigma^{-1}D\sigma$ for all $j_1,\ldots,j_n\in \mathbb{N}$, $j_1+\cdots+j_n\geq 1$. Thus, $\operatorname{Im}\sigma^{-1}D\sigma$ is an ideal generated by $x_1,\ldots,x_n$. Hence $\operatorname{Im}D$ is an ideal.
\end{proof}

\begin{prop}\label{prop2.7}
Let $R$ be a commutative ring and $\mathcal{A}$ an associative $R$-algebra. If $\mathfrak{r}(M)\subseteq M$ and $\mathfrak{r}(M)$ is an ideal of $\mathcal{A}$, then $M$ is a Mathieu-Zhao space of $R$.
\end{prop}
\begin{proof}
Since $\mathfrak{r}(M)$ is an ideal of $\mathcal{A}$, we have $\alpha^m\beta\in \mathfrak{r}(M)$ for any $\alpha\in \mathfrak{r}(M)$, $\beta\in \mathcal{A}$, $m\geq 1$. Since $\mathfrak{r}(M)\subseteq M$, we have $\alpha^m\beta\in M$ for any $\alpha\in \mathfrak{r}(M)$, $\beta\in \mathcal{A}$, $m\geq 1$. It follows from Proposition 2.1 in \cite{2} that $M$ is a Mathieu-Zhao space of $R$.
\end{proof}

\section{Conjecture \ref{conj1.1} for $\mathcal{E}$-derivations in dimension two}

\begin{thm}\label{thm3.1}
Let $\delta=I-\phi$ be an $\mathcal{E}$-derivation of $K[x]$. Then we have the following statements:

$(1)$ If $n=2r$ and $\phi(x_{2i-1})=\lambda_ix_{2i-1}+x_{2i}$, $\phi(x_{2i})=\lambda_ix_{2i}$, where $\lambda_i\in K$ for all $1\leq i\leq r$, then $\operatorname{Im}\delta$ is a Mathieu-Zhao space of $K[x]$.

$(2)$ If $n=2r+1$ and $\phi(x_{2i-1})=\lambda_ix_{2i-1}+x_{2i}$, $\phi(x_{2i})=\lambda_ix_{2i}$ and $\phi(x_{2r+1})=\lambda_{r+1}x_{2r+1}$, where $\lambda_i, \lambda_{r+1}\in K$ for all $1\leq i\leq r$, then $\operatorname{Im}\delta$ is a Mathieu-Zhao space of $K[x]$.

$(3)$ If $\phi(x_{2i-1})=\lambda_ix_{2i-1}+x_{2i}$, $\phi(x_{2i})=\lambda_ix_{2i}$ for all $1\leq i\leq t$ and $\phi(x_s)=\lambda_{s-t}x_s$ for all $2t+1\leq s\leq n$, where $1\leq t\leq \frac{n}{2}$, $t\in \mathbb{N}^*$, then $\operatorname{Im}\delta$ is a Mathieu-Zhao space of $K[x]$.
\end{thm}
\begin{proof}
$(1)$ Note that $\delta(x_{2r}^{i_{2r}})=(1-\lambda_r^{i_{2r}})x_{2r}^{i_{2r}}$. If $\lambda_r^{i_{2r}}\neq 1$, then $x_{2r}^{i_{2r}}\in \operatorname{Im}\delta$. If $\lambda_r^{i_{2r}}= 1$ for some $i_{2r}\in \mathbb{N}^*$, then $\delta(x_{2r-1}x_{2r}^{i_{2r}-1})=-\lambda_r^{i_{2r}-1}x_r^{i_2r}$. Thus, we have $x_r^{i_{2r}}\in \operatorname{Im}\delta$ for $i_{2r}\in \mathbb{N}^*$. Therefore, we have $x_{2r}^{i_{2r}}\in \operatorname{Im}\delta$ for all $i_{2r}\in \mathbb{N}^*$. Suppose that $x_{2r-1}^{l_{2r-1}}x_{2r}^{l_{2r}}\in \operatorname{Im}\delta$ for all $0\leq l_{2r-1}\leq i_{2r-1}-1$, $l_{2r}\in \mathbb{N}^*$. Then we have
\begin{eqnarray}
\nonumber
\delta(x_{2r-1}^{i_{2r-1}}x_{2r}^{i_{2r}})=x_{2r-1}^{i_{2r-1}}x_{2r}^{i_{2r}}-(\lambda_rx_{2r-1}+x_{2r})^{i_{2r-1}}(\lambda_rx_{2r})^{i_{2r}}~~~~~~~~~~~~~~~~~~~~~~~~~~~~~~\\
\nonumber
=(1-\lambda_r^{i_{2r-1}+i_{2r}})x_{2r-1}^{i_{2r-1}}x_{2r}^{i_{2r}}-\sum_{q_r=1}^{i_{2r-1}}\left(\begin{matrix}
i_{2r-1} \\
q_r
\end{matrix} \right)(\lambda_rx_{2r-1})^{i_{2r-1}-q_r}x_{2r}^{q_r}(\lambda_rx_{2r})^{i_{2r}}\\
\nonumber
 \end{eqnarray}
for $i_{2r}\in \mathbb{N}$. Note that $x_{2r-1}^{i_{2r-1}-q_r}x_{2r}^{q_r+i_{2r}}\in \operatorname{Im}\delta$ for all $1\leq q_r\leq i_{2r-1}$. If $\lambda_r^{i_{2r-1}+i_{2r}}\neq 1$, then $x_{2r-1}^{i_{2r-1}}x_{2r}^{i_{2r}}\in \operatorname{Im}\delta$. If $\lambda_r^{i_{2r-1}+i_{2r}}=1$, then we have
\begin{eqnarray}
\nonumber
\delta(x_{2r-1}^{i_{2r-1}+1}x_{2r}^{i_{2r}-1})=x_{2r-1}^{i_{2r-1}+1}x_{2r}^{i_{2r}-1}-(\lambda_rx_{2r-1}+x_{2r})^{i_{2r-1}+1}(\lambda_rx_{2r})^{i_{2r}-1}~~~~~~~\\
\nonumber
=-(i_{2r-1}+1)\lambda_r^{i_{2r-1}+i_{2r}-1}x_{2r-1}^{i_{2r-1}}x_{2r}^{i_{2r}}~~~~~~~~~~~~~~~~~~~~~~~~~~~~\\
\nonumber
-\sum_{q_r=2}^{i_{2r-1}+1}\left(\begin{matrix}
i_{2r-1}+1 \\
q_r
\end{matrix} \right)(\lambda_rx_{2r-1})^{i_{2r-1}-q_r+1}x_{2r}^{q_r}(\lambda_rx_{2r})^{i_{2r}-1}\\
\nonumber
 \end{eqnarray}
for $i_{2r}\in \mathbb{N}^*$. Since $x_{2r-1}^{i_{2r-1}-q_r+1}x_{2r}^{q_r+i_{2r}-1}\in \operatorname{Im}\delta$ for all $2\leq q_r\leq i_{2r-1}+1$, we have $x_{2r-1}^{i_{2r-1}}x_{2r}^{i_{2r}}\in \operatorname{Im}\delta$. Thus, we have $x_{2r-1}^{i_{2r-1}}x_{2r}^{i_{2r}}\in \operatorname{Im}\delta$ for all $i_{2r-1}\in \mathbb{N}$, $i_{2r}\in \mathbb{N}^*$. Since $\delta(x_2^{i_2}x_4^{i_4}\cdots x_{2r}^{i_{2r}})=(1-\lambda_1^{i_2}\lambda_2^{i_4}\cdots \lambda_r^{i_{2r}})x_2^{i_2}x_4^{i_4}\cdots x_{2r}^{i_{2r}}$, we have the following statements:\\
If $\lambda_1^{i_2}\lambda_2^{i_4}\cdots \lambda_r^{i_{2r}}\neq 1$ and $i_2+i_4+\cdots+i_{2r}\geq 1$, then $x_2^{i_2}x_4^{i_4}\cdots x_{2r}^{i_{2r}}\in \operatorname{Im}\delta$.\\
If $\lambda_1^{i_2}\lambda_2^{i_4}\cdots \lambda_r^{i_{2r}}= 1$ and $i_2+i_4+\cdots+i_{2r}\geq 1$, then we can assume that $i_2\geq 1$. Thus, we have
\begin{eqnarray}
\nonumber
\delta(x_1x_2^{i_2-1}x_4^{i_4}\cdots x_{2r}^{i_{2r}})=x_1x_2^{i_2-1}\cdots x_{2r}^{i_{2r}}-\lambda_1^{i_2-1}\lambda_2^{i_4}\cdots \lambda_r^{i_{2r}}(\lambda_1x_1+x_2)x_2^{i_2-1}\cdots x_{2r}^{i_{2r}}\\
\nonumber
=-\lambda_1^{i_2-1}\lambda_2^{i_4}\cdots \lambda_r^{i_{2r}}x_2^{i_2}x_4^{i_4}\cdots x_{2r}^{i_{2r}}~~~~~~~~~~~~~~~~~~~~~~~~~~~~~~~~~~~\\
\nonumber
 \end{eqnarray}
Hence we have $x_2^{i_2}x_4^{i_4}\cdots x_{2r}^{i_{2r}}\in \operatorname{Im}\delta$, whence $x_2^{i_2}x_4^{i_4}\cdots x_{2r}^{i_{2r}}\in \operatorname{Im}\delta$ for all $i_2+\cdots+i_{2r}\geq 1$.

$(1.1)$ Suppose that $x_{2j-1}^{l_{2j-1}}x_{2j}^{l_{2j}}\cdots x_{2k-1}^{l_{2k-1}}x_{2k}^{l_{2k}}\in \operatorname{Im}\delta$ for $l_{2j-1}\leq i_{2j-1}-1$ or $l_{2j-1}+l_{2j+1}\leq i_{2j-1}+i_{2j+1}-1$ or $\cdots$ or $l_{2j-1}+\cdots+l_{2k-1}\leq i_{2j-1}+\cdots+i_{2k-1}-1$, $l_{2j},\ldots,l_{2k}\in \mathbb{N}$ and $k$ is the maximal number such that $i_{2k}\neq 0$. We assume that if $i_{2m}=0$, then $i_{2m-1}=0$ for some $m\in\{j,j+1,\ldots,r\}$ and $1\leq j\leq r$. Then we have
\begin{eqnarray}
\nonumber
\delta(x_{2j-1}^{i_{2j-1}}\cdots x_{2k}^{i_{2k}})=x_{2j-1}^{i_{2j-1}}\cdots x_{2k}^{i_{2k}}-(\lambda_jx_{2j-1}+x_{2j})^{i_{2j-1}}\cdots(\lambda_kx_{2k})^{i_{2k}}\\
\nonumber
=(1-\lambda_j^{i_{2j-1}+i_{2j}}\cdots \lambda_k^{i_{2k-1}+i_{2k}} )x_{2j-1}^{i_{2j-1}}x_{2j}^{i_{2j}}\cdots x_{2k-1}^{i_{2k-1}}x_{2k}^{i_{2k}}+\operatorname{polynomial}~\operatorname{in}~\operatorname{Im}\delta
\end{eqnarray}
for $i_{2j-1},\ldots,i_{2k}\in \mathbb{N}$.\\
If $\lambda_j^{i_{2j-1}+i_{2j}}\cdots \lambda_k^{i_{2k-1}+i_{2k}}\neq 1$, then $x_{2j-1}^{i_{2j-1}}x_{2j}^{i_{2j}}\cdots x_{2k-1}^{i_{2k-1}}x_{2k}^{i_{2k}}\in \operatorname{Im}\delta$ for all $i_{2j-1},\ldots,\allowbreak i_{2k}\in \mathbb{N}$.\\
If $\lambda_j^{i_{2j-1}+i_{2j}}\cdots \lambda_k^{i_{2k-1}+i_{2k}}= 1$, then we have
\begin{eqnarray}
\nonumber
\delta(x_{2j-1}^{i_{2j-1}}\cdots x_{2k-1}^{i_{2k-1}+1}x_{2k}^{i_{2k}-1})=x_{2j-1}^{i_{2j-1}}\cdots x_{2k}^{i_{2k}}-(\lambda_jx_{2j-1}+x_{2j})^{i_{2j-1}}\cdots(\lambda_kx_{2k})^{i_{2k}-1}\\
\nonumber
=-(i_{2k-1}+1)\lambda_j^{i_{2j-1}+i_{2j}}\cdots \lambda_k^{i_{2k-1}+i_{2k}-1} x_{2j-1}^{i_{2j-1}}\cdots x_{2k}^{i_{2k}}+\operatorname{polynomial}~\operatorname{in}~\operatorname{Im}\delta
\end{eqnarray}
for $i_{2j-1},\ldots,i_{2k}\in \mathbb{N}$, $i_{2k}\geq 1$. Thus, we have $x_{2j-1}^{i_{2j-1}}\cdots x_{2k-1}^{i_{2k-1}}x_{2k}^{i_{2k}}\in \operatorname{Im}\delta$ for all $i_{2k}\geq 1$, $i_{2j-1},\ldots,i_{2k}\in \mathbb{N}$.

We have $x_{2j-2}^{i_{2j-2}}\cdots x_{2k-1}^{i_{2k-1}}x_{2k}^{i_{2k}}\in \operatorname{Im}\delta$ for all $i_{2k}\geq 1$, $i_{2j-2},\ldots,i_{2k}\in \mathbb{N}$ by the same arguments as above. Hence we have $x_1^{i_1}x_2^{i_2}\cdots x_{2k}^{i_{2k}}\in \operatorname{Im}\delta$ for all $i_{2k}\geq 1$ for some $k\in\{1,2,\ldots,r\}$, $i_1, i_2,\ldots,i_{2k-1}\in \mathbb{N}$.

$(1.2)$ Suppose that $x_{2j-1}^{l_{2j-1}}x_{2j}^{l_{2j}}\cdots x_{2\tilde{k}-2}^{l_{2\tilde{k}-2}}x_{2k-1}^{l_{2k-1}}\in \operatorname{Im}\delta$ for $l_{2j-1}\leq i_{2j-1}-1$ or $l_{2j-1}+l_{2j+1}\leq i_{2j-1}+i_{2j+1}-1$ or $\cdots$ or $l_{2j-1}+\cdots+l_{2k-1}\leq i_{2j-1}+\cdots+i_{2k-1}-1$, $l_{2j-1},\ldots,l_{2k-1}\in \mathbb{N}$ and $k$ is the maximal number such that $i_{2k-1}\neq 0$ and $i_{2\tilde{k}-2}\geq 1$ for some $\tilde{k}\in\{2,\ldots,r\}$, $\tilde{k}\leq k$ and $1\leq j\leq r$. Then we have
\begin{eqnarray}
\nonumber
\delta(x_{2j-1}^{i_{2j-1}}\cdots x_{2\tilde{k}-2}^{i_{2\tilde{k}-2}}x_{2k-1}^{i_{2k-1}})=x_{2j-1}^{i_{2j-1}}\cdots x_{2k-1}^{i_{2k-1}}-(\lambda_jx_{2j-1}+x_{2j})^{i_{2j-1}}\cdots(\lambda_kx_{2k-1}+x_{2k})^{i_{2k-1}}\\
\nonumber
=(1-\lambda_j^{i_{2j-1}+i_{2j}}\cdots \lambda_k^{i_{2k-1}} )x_{2j-1}^{i_{2j-1}}x_{2j}^{i_{2j}}\cdots x_{2\tilde{k}-2}^{i_{2\tilde{k}-2}}x_{2k-1}^{i_{2k-1}}+\operatorname{polynomial}~\operatorname{in}~\operatorname{Im}\delta~~~~~~
\end{eqnarray}
for $i_{2j-1},\ldots,i_{2\tilde{k}-2}, i_{2k-1}\in \mathbb{N}$.\\
If $\lambda_j^{i_{2j-1}+i_{2j}}\cdots \lambda_k^{i_{2k-1}}\neq 1$, then $x_{2j-1}^{i_{2j-1}}x_{2j}^{i_{2j}}\cdots x_{2\tilde{k}-2}^{i_{2\tilde{k}}-2}x_{2k-1}^{i_{2k-1}}\in \operatorname{Im}\delta$ for all $i_{2j-1},\ldots,\allowbreak i_{2\tilde{k}-2}, i_{2k-1}\in \mathbb{N}$.\\
If $\lambda_j^{i_{2j-1}+i_{2j}}\cdots \lambda_k^{i_{2k-1}}= 1$, then we have
\begin{eqnarray}
\nonumber
\delta(x_{2j-1}^{i_{2j-1}}\cdots x_{2\tilde{k}-3}^{i_{2\tilde{k}-3}+1} x_{2\tilde{k}-2}^{i_{2\tilde{k}-2}-1}x_{2k-1}^{i_{2k-1}})=x_{2j-1}^{i_{2j-1}}\cdots x_{2k-1}^{i_{2k-1}}-(\lambda_jx_{2j-1}+x_{2j})^{i_{2j-1}}\cdots\\
\nonumber
(\lambda_kx_{2k-1} +x_{2k})^{i_{2k-1}}\\
\nonumber
=-(i_{2\tilde{k}-3}+1)\lambda_j^{i_{2j-1}+i_{2j}}\cdots \lambda_k^{i_{2k-1}} x_{2j-1}^{i_{2j-1}}\cdots x_{2\tilde{k}-2}^{i_{2\tilde{k}-2}}x_{2k-1}^{i_{2k-1}}+Q_1(x)~~~~~~~~~~
\end{eqnarray}
for $i_{2j-1},\ldots,i_{2\tilde{k}-2}, i_{2k-1}\in \mathbb{N}$, $i_{2\tilde{k}-2}\geq 1$. Note that every monomial of $Q_1(x)$ is in $\operatorname{Im}\delta$ by the conclusion of $(1.1)$ and the induction hypothesis. Thus, we have
$x_{2j-1}^{i_{2j-1}}\cdots x_{2\tilde{k}-2}^{i_{2\tilde{k}-2}}x_{2k-1}^{i_{2k-1}}\in \operatorname{Im}\delta$ for all $i_{2\tilde{k}-2}\geq 1$, $i_{2j-1},\ldots,i_{2\tilde{k}-2},i_{2k-1}\in \mathbb{N}$.

We have $x_{2j-2}^{i_{2j-2}}\cdots x_{2\tilde{k}-2}^{i_{2\tilde{k}-2}}x_{2k-1}^{i_{2k-1}}\in \operatorname{Im}\delta$ for all $i_{2\tilde{k}-2}\geq 1$, $i_{2j-2},\ldots,i_{2\tilde{k}-2},i_{2k-1}\in \mathbb{N}$ by following the arguments of the former paragraph. Hence we have $x_1^{i_1}x_2^{i_2}\cdots \allowbreak x_{2\tilde{k}-2}^{i_{2\tilde{k}-2}} x_{2k-1}^{i_{2k-1}}\in \operatorname{Im}\delta$ for all $i_{2\tilde{k}-2}\geq 1$, $i_1,\ldots,i_{2\tilde{k}-2}, i_{2k-1}\in \mathbb{N}$.

$(1.3)$ We have that $x_{2j-1}^{i_{2j-1}}\cdots x_{2\tilde{k}-2}^{i_{2\tilde{k}-2}}x_{2k_1-1}^{i_{2k_1-1}}x_{2k-1}^{i_{2k-1}},~ x_{2j-2}^{i_{2j-2}}\cdots x_{2\tilde{k}-2}^{i_{2\tilde{k}-2}}x_{2k_1-1}^{i_{2k_1-1}}x_{2k-1}^{i_{2k-1}}\in \operatorname{Im}\delta$ for all $i_{2\tilde{k}-2}\geq 1$, $i_{2j-2}, i_{2j-1},\ldots,i_{2\tilde{k}-2}, i_{2k_1-1}, i_{2k-1}\in \mathbb{N}$, $\tilde{k}\leq k_1\leq k$ by following the arguments of $(1.2)$ and using the conclusions of $(1.1)$ and $(1.2)$. Thus, we have that $x_{2j-1}^{i_{2j-1}}\cdots x_{2\tilde{k}-2}^{i_{2\tilde{k}-2}}x_{2k_2-1}^{i_{2k_2-1}}\cdots x_{2k_1-1}^{i_{2k_1-1}}x_{2k-1}^{i_{2k-1}},~ x_{2j-2}^{i_{2j-2}}\cdots x_{2\tilde{k}-2}^{i_{2\tilde{k}-2}}x_{2k_2-1}^{i_{2k_2-1}}\cdots x_{2k_1-1}^{i_{2k_1-1}}\cdot\allowbreak x_{2k-1}^{i_{2k-1}}\in \operatorname{Im}\delta$ for all $i_{2\tilde{k}-2}\geq 1$, $i_{2j-2}, i_{2j-1},\ldots,i_{2\tilde{k}-2}, i_{2k_2-1},\ldots,i_{2k_1-1}, i_{2k-1}\in \mathbb{N}$, $\tilde{k}\leq k_2\leq\cdots\leq k_1\leq k$ by following the arguments of $(1.2)$ several times and using the former conclusions.

Combining the conclusion of $(1.1)$, we have that $x_1^{i_1}x_2^{i_2}\cdots x_{2r}^{i_{2r}}\in \operatorname{Im}\delta$ for all $i_1, i_2,\ldots,i_{2r}\in \mathbb{N}$ and $i_2+i_4+\cdots+i_{2r}\geq 1$. Thus, the ideal $I_1$ generated by $x_2, x_4,\ldots, x_{2r}$ is contained in $\operatorname{Im}\delta$.

Since
\begin{eqnarray}
\nonumber
\delta(x_1^{i_1}x_3^{i_3}\cdots x_{2r-1}^{i_{2r-1}})=(1-\lambda_1^{i_1}\lambda_2^{i_3}\cdots \lambda_r^{i_{2r-1}})x_1^{i_1}x_3^{i_3}\cdots x_{2r-1}^{i_{2r-1}} \operatorname{mod} I_1\\
\nonumber
=\hat{\delta}(x_1^{i_1}x_3^{i_3}\cdots x_{2r-1}^{i_{2r-1}})~~~~~~~~~~~~~~~~~~~~~~~~~~~~~~~~
\end{eqnarray}
for all $i_1, i_3,\ldots,i_{2r-1}\in \mathbb{N}$, where $\hat{\delta}=I-\hat{\phi}$ is an $\mathcal{E}$-derivation of $K[x_1,x_3,\ldots,x_{2r-1}]$ and $\hat{\phi}(x_{2j-1})=\lambda_jx_{2j-1}$ for all $1\leq j\leq r$, we have that $\operatorname{Im}\delta/I_1=\operatorname{Im}\hat{\delta}$. It follows from Lemma 3.2
and Corollary 3.3 in \cite{11} that $\operatorname{Im}\hat{\delta}$ is a Mathieu-Zhao space of $K[x_1,x_3,\ldots,x_{2r-1}]$. Then it follows from Proposition 2.7 in \cite{2} that $\operatorname{Im}\delta$ is a Mathieu-Zhao space of $K[x]$.\\

$(2)$ If $i_{2r+1}=0$, then we have that $x_1^{i_1}x_2^{i_2}\cdots x_{2r}^{i_{2r}}\in \operatorname{Im}\delta$ for all $i_2+\cdots+i_{2r}\geq 1$ by following the arguments of $(1)$.\\
If $i_{2r+1}\neq 0$, then we have that $x_1^{i_1}x_2^{i_2}\cdots x_{2r}^{i_{2r}}x_{2r+1}^{i_{2r+1}}\in \operatorname{Im}\delta$ for all $i_2+\cdots+i_{2r}\geq 1$ by following the arguments of $(1.2)$ and $(1.3)$. Thus, we have that the ideal $I_2$ generated by $x_2, x_4,\ldots,x_{2r}$ is contained in $\operatorname{Im}\delta$ and
\begin{eqnarray}
\nonumber
\delta(x_1^{i_1}x_3^{i_3}\cdots x_{2r+1}^{i_{2r+1}})=(1-\lambda_1^{i_1}\lambda_2^{i_3}\cdots \lambda_{r+1}^{i_{2r+1}})x_1^{i_1}x_3^{i_3}\cdots x_{2r+1}^{i_{2r+1}} \operatorname{mod} I_2\\
\nonumber
=\bar{\delta}(x_1^{i_1}x_3^{i_3}\cdots x_{2r+1}^{i_{2r+1}})~~~~~~~~~~~~~~~~~~~~~~~~~~~~~~~~
\end{eqnarray}
for all $i_1, i_3,\ldots,i_{2r+1}\in \mathbb{N}$, where $\bar{\delta}=I-\bar{\phi}$ is an $\mathcal{E}$-derivation of $K[x_1,x_3,\ldots,x_{2r+1}]$ and $\bar{\phi}(x_{2j-1})=\lambda_jx_{2j-1}$ for all $1\leq j\leq r+1$. Thus, we have $\operatorname{Im}\delta/I_2=\operatorname{Im}\bar{\delta}$. It follows from Lemma 3.2
and Corollary 3.3 in \cite{11} that $\operatorname{Im}\bar{\delta}$ is a Mathieu-Zhao space of $K[x_1,x_3,\ldots,x_{2r+1}]$. Then it follows from Proposition 2.7 in \cite{2} that $\operatorname{Im}\delta$ is a Mathieu-Zhao space of $K[x]$.\\

$(3)$ Following the arguments of $(2)$ by replacing $x_{2r+1}^{i_{2r+1}}$ with $x_{2t+1}^{i_{2t+1}}\cdots x_n^{i_n}$, we have that the ideal $I_3$ generated by $x_2, x_4,\ldots,x_{2t}$ is contained in $\operatorname{Im}\delta$ and
\begin{eqnarray}
\nonumber
\delta(x_1^{i_1}x_3^{i_3}\cdots x_{2t+1}^{i_{2t+1}}x_{2t+2}^{i_{2t+2}}\cdots x_n^{i_n})=(1-\lambda_1^{i_1}\lambda_2^{i_3}\cdots \lambda_{t+1}^{i_{2t+1}}\lambda_{t+2}^{i_{2t+2}}\cdots \lambda_{n-t}^{i_n})\cdot\\
\nonumber
x_1^{i_1}x_3^{i_3}\cdots x_{2t+1}^{i_{2t+1}}x_{2t+2}^{i_{2t+2}}\cdots x_n^{i_n} \operatorname{mod} I_3\\
\nonumber
=\tilde{\delta}(x_1^{i_1}x_3^{i_3}\cdots x_{2t+1}^{i_{2t+1}}x_{2t+2}^{i_{2t+2}}\cdots x_n^{i_n})~~~~~~~~
\end{eqnarray}
for all $i_1, i_3,\ldots,i_{2t+1}, i_{2t+2},\ldots,i_n\in \mathbb{N}$, where $\tilde{\delta}=I-\tilde{\phi}$ is an $\mathcal{E}$-derivation of $K[x_1,x_3,\ldots,x_{2t+1},x_{2t+2},\ldots,x_n]$ and $\tilde{\phi}(x_{2i-1})=\lambda_ix_{2i-1}$ for all $1\leq i\leq t$ and $\tilde{\phi}(x_s)=\lambda_{s-t}x_s$ for all $2t+1\leq s\leq n$. Thus, we have $\operatorname{Im}\delta/I_3=\operatorname{Im}\tilde{\delta}$. It follows from Lemma 3.2 and Corollary 3.3 in \cite{11} that $\operatorname{Im}\tilde{\delta}$ is a Mathieu-Zhao space of $K[x_1,x_3,\ldots,x_{2t+1},x_{2t+2},\ldots,x_n]$. Then it follows from Proposition 2.7 in \cite{2} that $\operatorname{Im}\delta$ is a Mathieu-Zhao space of $K[x]$.
\end{proof}

\begin{prop}\label{prop3.2}
Let $\delta=I-\phi$ be an $\mathcal{E}$-derivation of $K[x_1,x_2]$. If $\phi$ is a linear polynomial homomorphism of $K[x_1,x_2]$, then $\operatorname{Im}\delta$ is a Mathieu-Zhao space of $K[x_1,x_2]$.
\end{prop}
\begin{proof}
Since $\phi$ is a linear polynomial homomorphism, we have that
$$\left( \begin{matrix}
\phi(x_1)\\
\phi(x_2)
\end{matrix} \right)=A\left( \begin{matrix}
x_1 \\
x_2
\end{matrix} \right),$$
where $A\in M_2(K)$. Hence there exists $T\in \operatorname{GL}_2(K)$ such that
$$T^{-1}AT=\left( \begin{matrix}
\lambda_1 & 0\\
0 & \lambda_2
\end{matrix} \right)~\operatorname{or}~\left( \begin{matrix}
\lambda & 1\\
0 & \lambda
\end{matrix} \right),$$
where $\lambda_1\neq \lambda_2$. Let $(\sigma(x_1),\sigma(x_2))^t=T(x_1,x_2)^t$. Then we have $\sigma^{-1}\delta\sigma=I-\sigma^{-1}\phi\sigma$. It suffices to prove that $\operatorname{Im}(\sigma^{-1}\delta\sigma)$ is a Mathieu-Zhao space of $K[x_1,x_2]$. Let $\check{\delta}=\sigma^{-1}\delta\sigma=I-\check{\phi}$. Then $\check{\phi}(x_1)=\lambda_1x_1$, $\check{\phi}(x_2)=\lambda_2x_2$ or $\check{\phi}(x_1)=\lambda x_1+x_2$, $\check{\phi}(x_2)=\lambda x_2$.

$(1)$ If $\check{\phi}(x_1)=\lambda_1x_1$, $\check{\phi}(x_2)=\lambda_2x_2$, then it follows from Lemma 3.2 and Corollary 3.3 in \cite{11} that $\operatorname{Im}\check{\delta}$ is a Mathieu-Zhao space of $K[x_1,x_2]$.

$(2)$ If $\check{\phi}(x_1)=\lambda x_1+x_2$, $\check{\phi}(x_2)=\lambda x_2$, then it follows from Theorem \ref{thm3.1} $(1)$ that $\operatorname{Im}\check{\delta}$ is a Mathieu-Zhao space of $K[x_1,x_2]$. Then the conclusion follows.
\end{proof}

\begin{cor}\label{cor3.3}
Let $\delta=I-\phi$ be an $\mathcal{E}$-derivation of $K[x_1,x_2]$. If $\phi(x_1)=\lambda x_1+x_2$, $\phi(x_2)=\lambda x_2$, then $\operatorname{Im}\delta$ is an ideal or $\mathfrak{r}(\operatorname{Im}\delta)$ is an ideal of $K[x_1,x_2]$.
\end{cor}
\begin{proof}
$(1)$ If $\lambda$ is not a root of unity, then it follows from Theorem \ref{thm2.1} that $\operatorname{Im}\delta$ is an ideal of $K[x_1,x_2]$.

$(2)$ If $\lambda$ is a root of unity, then it follows from the proof of Theorem \ref{thm3.1} $(1)$ that $x_1^{i_1}x_2^{i_2}\in \operatorname{Im}\delta$ for all $i_1\in \mathbb{N}$, $i_2\in\mathbb{N}^*$ and $x_1^{i_1}\in \operatorname{Im}\delta$ for all $i_1\neq ds$, $d\in \mathbb{N}$, where $s$ is the least positive integer such that $\lambda^s=1$. That is, $x_1^{ds}\notin \operatorname{Im}\delta$ for all $d\in \mathbb{N}$. Next we prove that $\mathfrak{r}(\operatorname{Im}\delta)$ is an ideal generated by $x_2$. Clearly, the ideal generated by $x_2$ is contained in $\mathfrak{r}(\operatorname{Im}\delta)$. Let $G(x_1,x_2)=x_2G_1(x_1,x_2)+G_2(x_1)\in \mathfrak{r}(\operatorname{Im}\delta)$ and $G_2(x_1)\in K[x_1]$. We claim that $G_2(x_1)=0$. Otherwise, we have $G^m\in \operatorname{Im}\delta$ for all $m>>0$. Thus, we have $G_2^m\in \operatorname{Im}\delta$ for all $m>>0$. In particular, $G_2^{ds}\in \operatorname{Im}\delta$ for all $d>>0$. Suppose that $x_1^{\hat{t}}$ is the leading monomial of $G_2(x_1)$. Since $\operatorname{Im}\delta$ is a homogeneous $K$-subspace of $K[x_1,x_2]$, we have $x_1^{\hat{t}ds}\in \operatorname{Im}\delta$ for all $d>>0$, which is a contradiction. Thus, we have $G_2(x_1)=0$. Therefore, $G$ belongs to the ideal generated by $x_2$. Then the conclusion follows.
\end{proof}

\begin{prop}\label{prop3.4}
Let $\delta=I-\phi$ be an $\mathcal{E}$-derivation of $K[x_1,x_2]$. If $\phi$ is an affine polynomial homomorphism of $K[x_1,x_2]$, then $\operatorname{Im}\delta$ is a Mathieu-Zhao space of $K[x_1,x_2]$.
\end{prop}
\begin{proof}
Since $\phi$ is an affine polynomial homomorphism, we have that
$$\left( \begin{matrix}
\phi(x_1)\\
\phi(x_2)
\end{matrix} \right)=A\left( \begin{matrix}
x_1 \\
x_2
\end{matrix} \right)+\left( \begin{matrix}
c_1\\
c_2
\end{matrix} \right),$$
where $A\in M_2(K)$ and $(c_1,c_2)^t\in K^2$. Hence there exists $T\in \operatorname{GL}_2(K)$ such that
$$T^{-1}AT=\left( \begin{matrix}
\lambda_1 & 0\\
0 & \lambda_2
\end{matrix} \right)~\operatorname{or}~\left( \begin{matrix}
\lambda & 1\\
0 & \lambda
\end{matrix} \right),$$
where $\lambda_1\neq \lambda_2$. Let $(\sigma(x_1),\sigma(x_2))^t=T(x_1,x_2)^t$. Then we have $\sigma^{-1}\delta\sigma=I-\sigma^{-1}\phi\sigma$. It suffices to prove that $\operatorname{Im}(\sigma^{-1}\delta\sigma)$ is a Mathieu-Zhao space of $K[x_1,x_2]$. Let $\check{\delta}=\sigma^{-1}\delta\sigma=I-\check{\phi}$. Then $\check{\phi}(x_1)=\lambda_1x_1+\mu_1$, $\check{\phi}(x_2)=\lambda_2x_2+\mu_2$ or $\check{\phi}(x_1)=\lambda x_1+x_2+\mu_1$, $\check{\phi}(x_2)=\lambda x_2+\mu_2$, where $(\mu_1,\mu_2)^t=T^{-1}(c_1,c_2)^t$.

$(1)$ If $\lambda_1\neq 1$, $\lambda_2\neq 1$ and $\lambda\neq 1$, then it follows from Lemma \ref{lem2.2} that there exists $\check{\sigma}\in \operatorname{Aut}(K[x_1,x_2])$ such that $\check{\sigma}^{-1}\check{\delta}\check{\sigma}=I-\bar{\phi}$, where $\bar{\phi}$ is a linear polynomial homomorphism. Then it follows from Proposition \ref{prop3.2} that
$\operatorname{Im}(\check{\sigma}^{-1}\check{\delta}\check{\sigma})$ is a Mathieu-Zhao space of $K[x_1,x_2]$. Since $\check{\sigma}$ is a polynomial automorphism, we have that $\operatorname{Im}\check{\delta}$ is a Mathieu-Zhao space of $K[x_1,x_2]$.

$(2)$ If $\lambda_1=1$, then $\check{\phi}(x_1)=x_1+\mu_1$, $\check{\phi}(x_2)=\lambda_2x_2+\mu_2$. Thus, we have $\check{\delta}(x_1)=-\mu_1$. If $\mu_1\neq 0$, then $1\in \operatorname{Im}\check{\delta}$. It's easy to check that $\check{\delta}$ is locally finite. It follows from Proposition 1.4 in \cite{10} that $\operatorname{Im}\check{\delta}$ is a Mathieu-Zhao space of $K[x_1,x_2]$. If $\mu_1=0$, then $\check{\delta}(x_1^{i_1})=0$ for all $i_1\in \mathbb{N}$. Since $\lambda_2\neq \lambda_1$, there exists $\tau\in \operatorname{Aut}(K[x_1,x_2])$ such that $\tilde{\delta}:=\tau^{-1}\check{\delta}\tau=I-\tilde{\phi}$, where $\tilde{\phi}(x_1)=x_1$, $\tilde{\phi}(x_2)=\lambda_2x_2$. Then it follows from Proposition \ref{prop3.2} that $\operatorname{Im}\tilde{\delta}$ is a Mathieu-Zhao space of $K[x_1,x_2]$. Thus, $\operatorname{Im}\check{\delta}$ is a Mathieu-Zhao space of $K[x_1,x_2]$.

$(3)$ If $\lambda_2=1$, then we have that $\operatorname{Im}\check{\delta}$ is a Mathieu-Zhao space of $K[x_1,x_2]$ by following the arguments of Proposition 3.4 $(2)$.

$(4)$ If $\lambda=1$, then $\check{\phi}(x_1)=x_1+x_2+\mu_1$, $\check{\phi}(x_2)=x_2+\mu_2$. Thus, we have $\check{\delta}(x_2)=-\mu_2$. If $\mu_2\neq 0$, then $1\in \operatorname{Im}\check{\delta}$. Since $\check{\delta}$ is locally finite, it follows from Proposition 1.4 in \cite{10} that $\operatorname{Im}\check{\delta}$ is a Mathieu-Zhao space of $K[x_1,x_2]$. If $\mu_2=0$, then $\check{\delta}(x_2^{i_2})=0$ for all $i_2\in \mathbb{N}$. Thus, we have
$$\check{\delta}(x_1^{i_1}x_2^{i_2})=-(x_2+\mu_1)(\sum_{j=0}^{i_1-1}x_1^{i_1-j-1}(x_1+x_2+\mu_1)^j)x_2^{i_2}$$
for $i_1\in \mathbb{N}^*$, $i_2\in \mathbb{N}$. It's easy to check that $(x_2+\mu_1)x_1^{i_1}x_2^{i_2}\in \operatorname{Im}\check{\delta}$ for all $i_1,~i_2\in \mathbb{N}$. Since $1\notin \operatorname{Im}\check{\delta}$, we have that $\operatorname{Im}\check{\delta}$ is an ideal generated by $x_2+\mu_1$. Then the conclusion follows.
\end{proof}

\section{Conjecture \ref{conj1.1} for $\mathcal{E}$-derivations in dimension three}

\begin{lem}\label{lem4.1}
Let $\lambda_1, \lambda_2$ be elements in $K$. Then we have the following statements:

$(1)$ If one of $\lambda_1, \lambda_2$ is a root of unity and there exist $r_1, r_2\in \mathbb{N}^*$ such that $\lambda_1^{r_1}\lambda_2^{r_2}=1$, then the other is a root of unity.

$(2)$ If $\lambda_1^{r_1}\lambda_2^{r_2}=1$ and $\lambda_1^{\tilde{r}_1}\lambda_2^{\tilde{r}_2}=1$ for some $(\tilde{r}_1,\tilde{r}_2)\neq d(r_1,r_2)$ for any $d\in \mathbb{Q}^*$, $r_1,r_2,\tilde{r}_1,\tilde{r}_2\in \mathbb{N}^*$, then $\lambda_1, \lambda_2$ are both roots of unity.
\end{lem}
\begin{proof}
$(1)$ Without loss of generality, we can assume that $\lambda_1$ is a root of unity. Then there exists $s_1\in \mathbb{N}^*$ such that $\lambda_1^{s_1}=1$. Since $\lambda_1^{r_1}\lambda_2^{r_2}=1$, we have that $\lambda_2^{r_2}=\lambda_1^{-r_1}$. Hence we have $(\lambda_2^{r_2})^{s_1}=(\lambda_1^{s_1})^{-r_1}=1$. Since $r_2, s_1\in \mathbb{N}^*$, we have that $\lambda_2$ is a root of unity.

$(2)$ Since
\begin{equation}\label{eq4.1}
\lambda_1^{r_1}=\lambda_2^{-r_2}
\end{equation}
and
\begin{equation}\label{eq4.2}
(\lambda_1^{\tilde{r}_1}\lambda_2^{\tilde{r}_2})^{r_1}=1,
\end{equation}
we have $\lambda_2^{r_1\tilde{r}_2-\tilde{r}_1r_2}=1$ and $\lambda_2^{\tilde{r}_1r_2-r_1\tilde{r}_2}=1$ by substituting equation \eqref{eq4.1} to equation \eqref{eq4.2}. Since $r_1\tilde{r}_2\neq \tilde{r}_1r_2$, we have that $\lambda_2$ is a root of unity. It follows from Lemma \ref{lem4.1} $(1)$ that $\lambda_1$ is a root of unity.
\end{proof}

\begin{lem}\label{lem4.2}
Let $\delta=I-\phi$ be an $\mathcal{E}$-derivation of $K[x]$ and $\phi(x_i)=\lambda x_i+x_{i+1}$, $\phi(x_n)=\lambda x_n$ for $\lambda\in K$, $1\leq i\leq n-1$. If $\lambda$ is a root of unity, then $x_{n-1}^{i_{n-1}}x_n^{i_n}\in \operatorname{Im}\delta$ for all $i_{n-1}\in \mathbb{N}$, $i_n\in \mathbb{N}^*$ and $x_1^{i_1}x_2^{i_2}\cdots x_n^{i_n}\in \operatorname{Im}\delta$ in the case that $i_1+i_2+\cdots+i_n\neq ds$ for all $d\in \mathbb{N}$, where $s$ is the least positive integer such that $\lambda^s=1$.
\end{lem}
\begin{proof}
Since $\delta(x_n^{i_n})=(1-\lambda^{i_n})x_n^{i_n}$, we have $x_n^{i_n}\in \operatorname{Im}\delta$ for $i\neq ds$ for all $d\in \mathbb{N}$. Since $\delta(x_{n-1}x_n^{ds-1})=-\lambda^{ds-1}x_n^{ds}$, we have $x_n^{ds}\in \operatorname{Im}\delta$ for all $d\in \mathbb{N}^*$. Thus, we have $x_n^{i_n}\in \operatorname{Im}\delta$ for all $i_n\in \mathbb{N}^*$. Suppose that $x_{n-1}^{l_{n-1}}x_n^{i_n}\in \operatorname{Im}\delta$ for all $l_{n-1}\leq i_{n-1}-1$, $i_n\in \mathbb{N}^*$. Then
$$\delta(x_{n-1}^{i_{n-1}}x_n^{i_n})=(1-\lambda^{i_{n-1}+i_n})x_{n-1}^{i_{n-1}}x_n^{i_n}-\sum_{q_{n-1}=1}^{i_{n-1}}\left( \begin{matrix}
i_{n-1} \\
q_{n-1}
\end{matrix} \right)(\lambda x_{n-1})^{i_{n-1}-q_{n-1}}x_n^{q_{n-1}}(\lambda x_n)^{i_n}.$$
If $i_{n-1}+i_n\neq ds$ for all $d\in \mathbb{N}$, then we have $x_{n-1}^{i_{n-1}}x_n^{i_n}\in \operatorname{Im}\delta$.\\
If $i_{n-1}+i_n=ds$ for some $d\in \mathbb{N}^*$, then
$$\delta(x_{n-1}^{i_{n-1}+1}x_n^{i_n-1})=-(i_{n-1}+1)\lambda^{ds-1}x_{n-1}^{i_{n-1}}x_n^{i_n}-$$
$$\sum_{q_{n-1}=2}^{i_{n-1}+1}\left( \begin{matrix}
i_{n-1}+1 \\
q_{n-1}
\end{matrix} \right)(\lambda x_{n-1})^{i_{n-1}-q_{n-1}+1}x_n^{q_{n-1}}(\lambda x_n)^{i_n-1}.$$
Since $x_{n-1}^{i_{n-1}-q_{n-1}+1}x_n^{i_n+q_{n-1}-1}\in \operatorname{Im}\delta$ for all $2\leq q_{n-1}\leq i_{n-1}+1$, we have $x_{n-1}^{i_{n-1}}x_n^{i_n}\in \operatorname{Im}\delta$. Thus, we have $x_{n-1}^{i_{n-1}}x_n^{i_n}\in \operatorname{Im}\delta$ for all $i_n \in \mathbb{N}^*$.

Since
$$\delta(x_{n-1}^{i_{n-1}})=(1-\lambda^{i_{n-1}})x_{n-1}^{i_{n-1}}-\sum_{q_{n-1}=1}^{i_{n-1}}
\left(\begin{matrix}
i_{n-1} \\
q_{n-1}
\end{matrix} \right)(\lambda x_{n-1})^{i_{n-1}-q_{n-1}}x_n^{q_{n-1}}$$
and $x_{n-1}^{i_{n-1}-q_{n-1}}x_n^{q_{n-1}}\in \operatorname{Im}\delta$ for all $1\leq q_{n-1}\leq i_{n-1}$, we have $x_{n-1}^{i_{n-1}}\in \operatorname{Im}\delta$ for $i_{n-1}\neq ds$ for all $d\in \mathbb{N}$. Suppose that $x_k^{l_k}x_{k+1}^{l_{k+1}}\cdots x_{n-1}^{l_{n-1}}x_n^{i_n}\in \operatorname{Im}\delta$ for all $l_k\leq i_k-1$, $l_k+l_{k+1}\leq i_k+i_{k+1}-1$ or $\cdots$ or $l_k+\cdots+l_{n-1}\leq i_k+\cdots+i_{n-1}-1$ and $l_k+l_{k+1}+\cdots+l_{n-1}+i_n\neq ds$ for all $d\in \mathbb{N}$. Then
$$\delta(x_k^{i_k}\cdots x_{n-1}^{i_{n-1}}x_n^{i_n})=(1-\lambda^{i_k+\cdots+i_n})x_k^{i_k}\cdots x_{n-1}^{i_{n-1}}x_n^{i_n}+P(x_k,\ldots,x_n).$$
By induction hypothesis, we have $P(x_k,\ldots,x_n)\in \operatorname{Im}\delta$ if $i_k+\cdots+i_n\neq ds$ for all $d\in \mathbb{N}$. Thus, we have $x_k^{i_k}\cdots x_{n-1}^{i_{n-1}}x_n^{i_n}\in \operatorname{Im}\delta$ for $i_k+\cdots+i_n\neq ds$ for all $d\in \mathbb{N}$.
\end{proof}

\begin{prop}\label{prop4.3}
Let $\delta=I-\phi$ be an $\mathcal{E}$-derivation of $K[x_1,x_2,x_3]$ and $\phi(x_1)=\lambda_1 x_1+x_2$, $\phi(x_2)=\lambda_1 x_2$ and $\phi(x_3)=\lambda_2x_3$ for $\lambda_1, \lambda_2\in K$. Then we have the following statements:

$(1)$ If $\lambda_1$ is a root of unity and $\lambda_2$ is not a root of unity, then $\mathfrak{r}(\operatorname{Im}\delta)$ is an ideal generated by $x_2, x_3$.

$(2)$ If $\lambda_1$ is not a root of unity and $\lambda_2$ is a root of unity, then $\mathfrak{r}(\operatorname{Im}\delta)$ is an ideal generated by $x_1, x_2$.
\end{prop}
\begin{proof}
It follows from Theorem \ref{thm3.1} $(3)$ that the ideal $(x_2)\subseteq \operatorname{Im}\delta$.

$(1)$ Since $\delta(x_3^{i_3})=(1-\lambda_2^{i_3})x_3^{i_3}$ and $\lambda_2^{i_3}\neq 1$ for any $i_3\in \mathbb{N}^*$, we have $x_3^{i_3}\in \operatorname{Im}\delta$ for all $i_3\in \mathbb{N}^*$. Suppose that $x_1^{l_1}x_3^{i_3}\in \operatorname{Im}\delta$ for all $0\leq l_1\leq i_1-1$, $i_3\geq 1$. Then we have
\begin{eqnarray}
\nonumber
\delta(x_1^{i_1}x_3^{i_3})=x_1^{i_1}x_3^{i_3}-(\lambda_1x_1+x_2)^{i_1}(\lambda_2x_3)^{i_3}\\
\nonumber
=(1-\lambda_1^{i_1}\lambda_2^{i_3})x_1^{i_1}x_3^{i_3}-\sum_{q_1=1}^{i_1}\left(\begin{matrix}
i_1 \\
q_1
\end{matrix} \right)(\lambda_1x_1)^{i_1-q_1}x_2^{q_1}(\lambda_2x_3)^{i_3}
\end{eqnarray}
for $i_1\in \mathbb{N}$, $i_3\in \mathbb{N}^*$. Since $(x_2)\subseteq \operatorname{Im}\delta$, we have $x_1^{i_1-q_1}x_2^{q_1}x_3^{i_3}\in \operatorname{Im}\delta$ for all $1\leq q_1\leq i_1$. It follows from Lemma \ref{lem4.1} $(1)$ that $\lambda_1^{i_1}\lambda_3^{i_3}\neq 1$ for all $i_3\in \mathbb{N}^*$. Thus, we have $x_1^{i_1}x_3^{i_3}\in \operatorname{Im}\delta$ for all $i_3\in \mathbb{N}^*$. Therefore, it follows from Lemma \ref{lem4.2} that $\operatorname{Im}\delta$ is a $K$-vector space generated by monomials $x_1^{i_1}x_2^{i_2}x_3^{i_3}$ for $i_1, i_2, i_3\in \mathbb{N}$, $i_2+i_3\geq 1$ and $x_1^{i_1}$ for $i_1\neq d_1s_1$ for all $d_1\in \mathbb{N}$, where $s_1$ is the least positive integer such that $\lambda_1^{s_1}=1$.

We claim that $g^m\in \operatorname{Im}\delta$ for any $g\in K[x_1,x_2,x_3]$, $m>>0$ iff $g(x_1,0,0)=0$. If $g(x_1,0,0)=0$, then it's easy to see that $g^m\in \operatorname{Im}\delta$ for all $m\geq 1$. Conversely, let $g(x_1,x_2,x_3)=g_1(x_1,x_2,x_3)+g_2(x_1)$, where $g_1\in K[x_1,x_2,x_3]$, $g_2\in K[x_1]$ and $g_1(x_1,0,0)=0$. If $g_2(x_1)\neq 0$, then we can assume that $x_1^{\bar{t}}$ be the leading term in $g_2(x_1)$. Since $g^m\in \operatorname{Im}\delta$ for all $m>>0$ and $\operatorname{Im}\delta$ is homogeneous, we have $x_1^{m\bar{t}}\in \operatorname{Im}\delta$ for all $m>>0$, which is a contradiction. Thus, we have $g_2(x_1)=0$. That is, $g(x_1,0,0)=0$. Therefore, $\mathfrak{r}(\operatorname{Im}\delta)$ is an ideal generated by $x_2, x_3$.

$(2)$ Since $\delta(x_3^{i_3})=(1-\lambda_2^{i_3})x_3^{i_3}$, we have that $x_3^{i_3}\in \operatorname{Im}\delta$ for $i_3\neq d_2s_2$ for all $d_2\in \mathbb{N}$, where $s_2$ is the least positive integer such that $\lambda_2^{s_2}=1$. Note that $\delta(x_1x_3^{i_3})=(1-\lambda_1\lambda_2^{i_3})x_1x_3^{i_3}-\lambda_2^{i_3}x_2x_3^{i_3}$ and the ideal $(x_2)\subseteq \operatorname{Im}\delta$. It follows from Lemma \ref{lem4.1} $(1)$ that $\lambda_1\lambda_2^{i_3}\neq 1$ for all $i_3\in \mathbb{N}$. Thus, we have $x_1x_3^{i_3}\in \operatorname{Im}\delta$ for all $i_3\in \mathbb{N}$. Suppose that $x_1^{l_1}x_3^{i_3}\in \operatorname{Im}\delta$ for all $1\leq l_1\leq i_1-1$, $i_3\in \mathbb{N}$. Then we have
\begin{eqnarray}
\nonumber
\delta(x_1^{i_1}x_3^{i_3})=x_1^{i_1}x_3^{i_3}-(\lambda_1x_1+x_2)^{i_1}(\lambda_2x_3)^{i_3}\\
\nonumber
=(1-\lambda_1^{i_1}\lambda_2^{i_3})x_1^{i_1}x_3^{i_3}-\sum_{q_1=1}^{i_1}\left(\begin{matrix}
i_1 \\
q_1
\end{matrix} \right)(\lambda_1x_1)^{i_1-q_1}x_2^{q_1}(\lambda_2x_3)^{i_3}
\end{eqnarray}
for $i_1\in \mathbb{N}^*$, $i_3\in \mathbb{N}$. Since $(x_2)\subseteq \operatorname{Im}\delta$, we have $x_1^{i_1-q_1}x_2^{q_1}x_3^{i_3}\in \operatorname{Im}\delta$ for all $1\leq q_1\leq i_1$. It follows from Lemma \ref{lem4.1} $(1)$ that $\lambda_1^{i_1}\lambda_3^{i_3}\neq 1$. Thus, we have $x_1^{i_1}x_3^{i_3}\in \operatorname{Im}\delta$ for all $i_1\in \mathbb{N}^*$. Therefore, $\operatorname{Im}\delta$ is a $K$-vector space generated by monomials $x_1^{i_1}x_2^{i_2}x_3^{i_3}$ for $i_1, i_2, i_3\in \mathbb{N}$, $i_1+i_2\geq 1$ and $x_3^{i_3}$ for $i_3\neq d_2s_2$ for all $d_2\in \mathbb{N}$.

We claim that $G^m\in \operatorname{Im}\delta$ for any $G\in K[x_1,x_2,x_3]$, $m>>0$ iff $G(0,0,x_3)=0$. If $G(0,0,x_3)=0$, then it's easy to see that $G^m\in \operatorname{Im}\delta$ for all $m\geq 1$. Conversely, let $G(x_1,x_2,x_3)=G_1(x_1,x_2,x_3)+G_2(x_3)$, where $G_1\in K[x_1,x_2,x_3]$, $G_2\in K[x_3]$ and $G_1(0,0,x_3)=0$. If $G_2(x_3)\neq 0$, then we can assume that $x_3^{\tilde{t}}$ be the leading term in $G_2(x_3)$. Since $G^m\in \operatorname{Im}\delta$ for all $m>>0$ and $\operatorname{Im}\delta$ is homogeneous, we have $x_3^{m\tilde{t}}\in \operatorname{Im}\delta$ for all $m>>0$, which is a contradiction. Thus, we have $G_2(x_3)=0$. That is, $G(0,0,x_3)=0$. Therefore, $\mathfrak{r}(\operatorname{Im}\delta)$ is an ideal generated by $x_1, x_2$.
\end{proof}

\begin{prop}\label{prop4.4}
Let $\delta=I-\phi$ be an $\mathcal{E}$-derivation of $K[x_1,x_2,x_3]$ and $\phi(x_i)=\lambda x_i+x_{i+1}$ and $\phi(x_3)=\lambda x_3$ for $\lambda\in K$ and $i=1, 2$. Then we have the following statements:

$(1)$ If $\lambda=1$, then $x_1^{i_1}x_2^{2k+1}x_3^{i_3}\in \operatorname{Im}\delta$ for $i_3\geq i_1\geq 0$, $x_1^{i_1}x_2^{2k}x_3^{i_3}\in \operatorname{Im}\delta$ for $i_3\geq i_1+1\geq 1$ and $x_1^{i_1}x_2^{2k+1}x_3^{i_3}\notin \operatorname{Im}\delta$ for $i_3<i_1$, $x_1^{i_1}x_2^{2k}x_3^{i_3}\notin \operatorname{Im}\delta$ for $i_3\leq i_1$. In particular, $x_1^{i_1}x_2^{i_2}x_3^{i_3}\in \mathfrak{r}(\operatorname{Im}\delta)$ for all $i_3\geq i_1+1$, $i_1,~i_2,~i_3,~ k\in\mathbb{N}$.

$(2)$ If $\lambda$ is a root of unity and $\lambda\neq 1$, then $x_1^{i_1}x_2^{i_2}x_3^{i_3}\in \operatorname{Im}\delta$ in the case that $i_1+i_2+i_3\neq ds$ for all $d\in \mathbb{N}$ and $x_1^{i_1}x_2^{i_2}x_3^{i_3}\in \operatorname{Im}\delta$ for $i_3\geq i_1+1\geq 1$ if $i_1+i_2+i_3=ds$ for some $d\in \mathbb{N}^*$, $x_1^{i_1}x_2^{i_2}x_3^{i_3}\notin \operatorname{Im}\delta$ for $i_3\leq i_1$ if $i_1+i_2+i_3=ds$ for some $d\in \mathbb{N}$, where $i_1,~i_2,~i_3\in \mathbb{N}$ and $s$ is the least positive integer such that $\lambda^s=1$. In particular, $x_1^{i_1}x_2^{i_2}x_3^{i_3}\in \mathfrak{r}(\operatorname{Im}\delta)$ for all $i_3\geq i_1+1$, $i_1,~i_2,~i_3\in\mathbb{N}$.
\end{prop}
\begin{proof}
$(1)$ It follows from Lemma \ref{lem4.2} that $x_2^{i_2}x_3^{i_3}\in \operatorname{Im}\delta$ for all $i_2\in \mathbb{N}$, $i_3\in \mathbb{N}^*$. Suppose that $x_1^{l_1}x_2^{l_2}x_3^{i_3}\in \operatorname{Im}\delta$ for $l_1\leq i_1-1$, $i_3\geq l_1+1$ or $l_1+l_2\leq i_1+i_2-1$ and $i_3\geq l_1+1$. Then we have
\begin{eqnarray}
\nonumber
\delta(x_1^{i_1}x_2^{i_2+1}x_3^{i_3})=x_1^{i_1}x_2^{i_2+1}x_3^{i_3}-(x_1+x_2)^{i_1}(x_2+x_3)^{i_2+1}x_3^{i_3}\\
\nonumber
=-(i_2+1)x_1^{i_1}x_2^{i_2}x_3^{i_3+1}+\operatorname{polynomial}~\operatorname{in}~ \operatorname{Im}\delta,
\end{eqnarray}
where $i_1,~i_2,~i_3\in \mathbb{N}$. Thus, we have $x_1^{i_1}x_2^{i_2}x_3^{i_3+1}\in \operatorname{Im}\delta$, where $i_3+1\geq (i_1-1)+1+1=i_1+1$. Hence we have $x_1^{i_1}x_2^{i_2}x_3^{i_3}\in \operatorname{Im}\delta$ for all $i_3\geq i_1+1$, $i_1,~i_2,~i_3\in\mathbb{N}$. Since $\delta(x_1)=-x_2$ and $\delta(x_1^2x_3^{i_3})=-(2x_1x_2+x_2^2)x_3^{i_3}$, we have $x_2,~x_1x_2x_3^{i_3}\in \operatorname{Im}\delta$ for all $i_3\in \mathbb{N}^*$. Suppose that $x_1^{l_1}x_2^{2l_2-1}x_3^{i_3}\in \operatorname{Im}\delta$ for $l_1\leq i_1-1$, $i_3\geq l_1$ or $l_1+l_2\leq i_1+k-1$, $i_3\geq l_1$. Then we have
\begin{eqnarray}
\nonumber
\delta(x_1^{i_1}x_2^{2k}x_3^{i_3})=x_1^{i_1}x_2^{2k}x_3^{i_3}-(x_1+x_2)^{i_1}(x_2+x_3)^{2k}x_3^{i_3}\\
\nonumber
=-(2k)x_1^{i_1}x_2^{2k-1}x_3^{i_3+1}+P(x_1,x_2,x_3),
\end{eqnarray}
where every monomial of $P(x_1,x_2,x_3)$ belongs to $\operatorname{Im}\delta$ and $i_1,~i_3\in \mathbb{N}$, $k\in \mathbb{N}^*$. Thus, we have $x_1^{i_1}x_2^{2k-1}x_3^{i_3+1}\in \operatorname{Im}\delta$, where $i_3+1\geq (i_1-1)+1=i_1$ and $i_1,~i_3\in \mathbb{N}$, $k\in \mathbb{N}^*$. Hence we have $x_1^{i_1}x_2^{2k-1}x_3^{i_3}\in \operatorname{Im}\delta$ for all $i_3\geq i_1$, $i_1,~i_3\in \mathbb{N}$, $k\in \mathbb{N}^*$.

If $0\leq i_3\leq i_1-1$, then we have
\begin{eqnarray}\label{eq4.3}
\delta(x_1^{i_1}x_2^{2k+1}x_3^{i_3})=x_1^{i_1}x_2^{2k+1}x_3^{i_3}-(x_1+x_2)^{i_1}(x_2+x_3)^{2k+1}x_3^{i_3}\\
\nonumber
=x_1^{i_1}x_2^{2k+1}x_3^{i_3}-\sum_{q_1=0}^{i_1}\sum_{q_2=0}^{i_1-i_3-q_1}\left(\begin{matrix}
i_1 \\
q_1
\end{matrix} \right)\left(\begin{matrix}
2k+1 \\
q_2
\end{matrix} \right)x_1^{i_1-q_1}x_2^{2k-q_2+q_1+1}x_3^{i_3+q_2}\\
\nonumber
+\operatorname{polynomial}~ \operatorname{in}~\operatorname{Im}\delta
\end{eqnarray}
If $0\leq i_3\leq i_1-2$, then we have
\begin{eqnarray}\label{eq4.4}
\delta(x_1^{i_1}x_2^{2k}x_3^{i_3})=x_1^{i_1}x_2^{2k}x_3^{i_3}-(x_1+x_2)^{i_1}(x_2+x_3)^{2k}x_3^{i_3}\\
\nonumber
=x_1^{i_1}x_2^{2k}x_3^{i_3}-\sum_{q_1=0}^{i_1}\sum_{q_2=0}^{i_1-i_3-q_1}\left(\begin{matrix}
i_1 \\
q_1
\end{matrix} \right)\left(\begin{matrix}
2k \\
q_2
\end{matrix} \right)x_1^{i_1-q_1}x_2^{2k-q_2+q_1}x_3^{i_3+q_2}\\
\nonumber
+\operatorname{polynomial}~ \operatorname{in}~\operatorname{Im}\delta
\end{eqnarray}
for all $i_1,~i_3,~k \in \mathbb{N}$. It follows from equation \eqref{eq4.3} that at least two distinct monomials in $\delta(x_1^{i_1}x_2^{2k+1}x_3^{i_3})$ by module
$\operatorname{Im}\delta$ for all $0\leq i_3\leq i_1-1$ and
$S_{k,i_1,i_3}^{(1)}:=\{x_1^{i_1}x_2^{2k-i_1+i_3+1}x_3^{i_1},~x_1^{i_1-1}x_2^{2k-i_1+i_3+3}x_3^{i_1-1}\}\subseteq \{\operatorname{monomials}~\operatorname{in}~\delta(x_1^{i_1}x_2^{2k+1}x_3^{i_3})\} \operatorname{mod}\allowbreak \operatorname{Im}\delta$ if $i_1-i_3-1$ is even and $S_{k,i_1,i_3}^{(2)}:=\{x_1^{i_1}x_2^{2k-i_1+i_3+2}x_3^{i_1-1},~x_1^{i_1-1}x_2^{2k-i_1+i_3+4}x_3^{i_1-2}\}\subseteq \{\operatorname{monomials}~\operatorname{in}~\delta(x_1^{i_1}x_2^{2k+1}x_3^{i_3})\} \operatorname{mod}~\operatorname{Im}\delta$ if $i_1-i_3-1$ is odd for all $0\leq i_3\leq i_1-1$. It follows from equation \eqref{eq4.4} that at least two distinct monomials in $\delta(x_1^{i_1}x_2^{2k}x_3^{i_3})$ by module
$\operatorname{Im}\delta$ for all $0\leq i_3\leq i_1-2$ and
$S_{k,i_1,i_3}^{(3)}:=\{x_1^{i_1}x_2^{2k-i_1+i_3+2}x_3^{i_1-2},\allowbreak x_1^{i_1}x_2^{2k-i_1+i_3+1}x_3^{i_1-1}\}\subseteq \{\operatorname{monomials}~\operatorname{in}~\delta(x_1^{i_1}x_2^{2k}x_3^{i_3})\} \operatorname{mod}\allowbreak \operatorname{Im}\delta$ if $i_1-i_3-1$ is even and $S_{k,i_1,i_3}^{(4)}:=\{x_1^{i_1}x_2^{2k-i_1+i_3+1}x_3^{i_1-1},~x_1^{i_1}x_2^{2k-i_1+i_3}x_3^{i_1}\}\subseteq \{\operatorname{monomials}~\operatorname{in}~\delta(x_1^{i_1}x_2^{2k}x_3^{i_3})\} \allowbreak \operatorname{mod}~\operatorname{Im}\delta$ if $i_1-i_3-1$ is odd for all $0\leq i_3\leq i_1-2$. It's easy to check that $S_{k,i_1,i_3}^{(j_1)}\neq S_{\tilde{k},\tilde{i}_1,\tilde{i}_3}^{(j_2)}$ for $j_1\neq j_2$, $1\leq j_1,~j_2\leq 4$, $k,~\tilde{k},~i_1,~\tilde{i}_1,~i_3,~\tilde{i}_3 \in \mathbb{N}$ and $S_{k,i_1,i_3}^{(j)}\neq S_{\tilde{k},\tilde{i}_1,\tilde{i}_3}^{(j)}$ for $k\neq \tilde{k}$ and $S_{k,i_1,i_3}^{(j)}\neq S_{k,\tilde{i}_1,\tilde{i}_3}^{(j)}$ for $i_1\neq \tilde{i}_1$ or $i_3\neq \tilde{i}_3$ for all $1\leq j\leq 4$. Thus, any linear combination of $\delta(x_1^{i_1}x_2^{2k+1}x_3^{i_3})$ for all $0\leq i_3\leq i_1-1$ and $\delta(x_1^{i_1}x_2^{2k}x_3^{i_3})$ for all $0\leq i_3\leq i_1-2$ has at least two distinct monomials by module $\operatorname{Im}\delta$. Hence $x_1^{i_1}x_2^{2k+1}x_3^{i_3}\notin \operatorname{Im}\delta$ for all $i_3<i_1$
and $x_1^{i_1}x_2^{2k}x_3^{i_3}\notin \operatorname{Im}\delta$ for all $i_3\leq i_1$, $i_1,~i_2,~i_3\in \mathbb{N}$.

$(2)$ Since $\lambda \neq 1$, we have $s\geq 2$. It follows from Lemma \ref{lem4.2} that $x_1^{i_1}x_2^{i_2}x_3^{i_3}\in \operatorname{Im}\delta$ if $i_1+i_2+i_3\neq ds$ for all $d\in \mathbb{N}$ and $x_2^{i_2}x_3^{i_3}\in \operatorname{Im}\delta$ for all $i_2\in \mathbb{N}$, $i_3\in \mathbb{N}^*$. Suppose that $x_1^{l_1}x_2^{ds-l_1-i_3}x_3^{i_3}\in \operatorname{Im}\delta$ for all $0\leq l_1\leq i_1-1$, $i_3\geq l_1+1$. Then we have
\begin{eqnarray}
\nonumber
\delta(x_1^{i_1}x_2x_3^{ds-i_1-1})=x_1^{i_1}x_2x_3^{ds-i_1-1}-(\lambda x_1+x_2)^{i_1}(\lambda x_2+x_3)(\lambda x_3)^{ds-i_1-1}\\
\nonumber
=-\lambda^{ds-1}x_1^{i_1}x_3^{ds-i_1}+\operatorname{polynomial}~\operatorname{in}~ \operatorname{Im}\delta,
\end{eqnarray}
where $ds-i_1-1=i_3\geq (i_1-1)+1=i_1$. Thus, we have $x_1^{i_1}x_3^{ds-i_1}\in \operatorname{Im}\delta$ for all $ds-i_1\geq i_1+1$. Suppose that $x_1^{l_1}x_2^{ds-l_1-i_3}x_3^{i_3}\in \operatorname{Im}\delta$ for all $0\leq l_1\leq i_1-1$, $i_3\geq l_1+1$ or $l_1=i_1$ and $i_3\geq i_1+\tilde{r}$, $\tilde{r}\geq 2$. Then we have
\begin{eqnarray}
\nonumber
\delta(x_1^{i_1}x_2^{ds-2i_1-\tilde{r}+2}x_3^{i_1+\tilde{r}-2})~~~~~~~~~~~~~~~~~~~~~~~~~~~~~~~~~~~~~~~~~~~~~~~~~~~~~~~~~~\\
\nonumber
=x_1^{i_1}x_2^{ds-2i_1-\tilde{r}+2}x_3^{i_1+\tilde{r}-2}-(\lambda x_1+x_2)^{i_1}(\lambda x_2+x_3)^{ds-2i_1-\tilde{r}+2}(\lambda x_3)^{i_1+\tilde{r}-2}\\
\nonumber
=-\lambda^{ds-1}(ds-2i_1-\tilde{r}+2)x_1^{i_1}x_2^{ds-2i_1-\tilde{r}+1}x_3^{i_1+\tilde{r}-1}+\operatorname{polynomial}~\operatorname{in}~ \operatorname{Im}\delta.
\end{eqnarray}
Thus, we have $x_1^{i_1}x_2^{ds-2i_1-\tilde{r}+1}x_3^{i_1+\tilde{r}-1}\in \operatorname{Im}\delta$ for all $\tilde{r}\geq 2$, $i_1\in \mathbb{N}$, $d\in \mathbb{N}^*$. Hence we have $x_1^{i_1}x_2^{i_2}x_3^{i_3}\in \operatorname{Im}\delta$ for all $i_3\geq i_1+1$ if $i_1+i_2+i_3=ds$ for some $d\in \mathbb{N}^*$.

If $0\leq i_3\leq i_1-1$, then we have
\begin{eqnarray}\label{eq4.5}
\delta(x_1^{i_1}x_2^{ds-i_1-i_3}x_3^{i_3})~~~~~~~~~~~~~~~~~~~~~~~~~~~~~~~~~~~~~~~~~~~~~~~~~~~~~~~~~~~~~~~~~~~~\\
\nonumber
=x_1^{i_1}x_2^{ds-i_1-i_3}x_3^{i_3}-(\lambda x_1+x_2)^{i_1}(\lambda x_2+x_3)^{ds-i_1-i_3}(\lambda x_3)^{i_3}\\
\nonumber
=x_1^{i_1}x_2^{ds-i_1-i_3}x_3^{i_3}-\sum_{q_1=0}^{i_1}\sum_{q_2=0}^{i_1-i_3-q_1}\left(\begin{matrix}
i_1 \\
q_1
\end{matrix} \right)\left(\begin{matrix}
ds-i_1-i_3 \\
q_2
\end{matrix} \right)x_1^{i_1-q_1}x_2^{ds-i_1-i_3-q_2+q_1}x_3^{i_3+q_2}\\
\nonumber
+\operatorname{polynomial}~\operatorname{in}~ \operatorname{Im}\delta
\end{eqnarray}
for all $i_1,~i_3\in \mathbb{N}$. It follows from equation \eqref{eq4.5} that at least two distinct monomials in $\delta(x_1^{i_1}x_2^{ds-i_1-i_3}x_3^{i_3})$ by module $\operatorname{Im}\delta$ for all $0\leq i_3\leq i_1-1$ and $S_{d,i_1,i_3}:=\{x_1^{i_1}x_2^{ds-2i_1}x_3^{i_1},~x_1^{i_1-1}x_2^{ds-2i_1+2}x_3^{i_1-1}\}\subseteq \{\operatorname{monomial} \operatorname{in} \delta(x_1^{i_1}x_2^{ds-i_1-i_3}x_3^{i_3})\}$ $\operatorname{mod} \operatorname{Im}\delta$
for all $0\leq i_3\leq i_1-1$. It's easy to check that $S_{d,i_1,i_3}\neq S_{\tilde{d},\tilde{i}_1,\tilde{i}_3}$ for $d\neq \tilde{d}$ and $S_{d,i_1,i_3}\neq S_{d,\tilde{i}_1,\tilde{i}_3}$ for $i_1\neq \tilde{i}_1$ or $i_3\neq \tilde{i}_3$. Thus, any linear combination of $\delta(x_1^{i_1}x_2^{ds-i_1-i_3}x_3^{i_3})$ for all $0\leq i_3\leq i_1-1$ has at least two distinct monomials by module $\operatorname{Im}\delta$. Hence $x_1^{i_1}x_2^{i_2}x_3^{i_3}\notin \operatorname{Im}\delta$ for all $i_3\leq i_1$ and $i_1+i_2+i_3=ds$ for some $d\in \mathbb{N}$.
\end{proof}

\begin{conj}\label{conj4.5}
Let $\delta=I-\phi$ be an $\mathcal{E}$-derivation of $K[x_1,x_2,x_3]$ and $\phi(x_i)=\lambda x_i+x_{i+1}$ and $\phi(x_3)=\lambda x_3$ for $\lambda\in K$, $i=1, 2$. Then $\mathfrak{r}(\operatorname{Im}\delta)$ is a $K$-vector space generated by the monomials $x_1^{i_1}x_2^{i_2}x_3^{i_3}$ for all $i_3\geq i_1+1$, $i_1,~i_2,~i_3 \in \mathbb{N}$.
\end{conj}

\begin{rem}\label{rem4.6}
It follows from Proposition \ref{prop4.4} that the $K$-vector space $V$ generated by the monomials $x_1^{i_1}x_2^{i_2}x_3^{i_3}$ for all $i_3\geq i_1+1$, $i_1,~i_2,~i_3 \in \mathbb{N}$ is contained in $\mathfrak{r}(\operatorname{Im}\delta)$. If $\mathfrak{r}(\operatorname{Im}\delta)=V$, then $\operatorname{Im}\delta$ is a Mathieu-Zhao space of $K[x_1,x_2,x_3]$ because for any $f\in \mathfrak{r}(\operatorname{Im}\delta)$, $h\in K[x_1,x_2,x_3]$, we have $hf^m\in \operatorname{Im}\delta$ for all $m\geq N$, where $N=\deg h+1$.

If $\lambda=1$, then $\delta$ is locally nilpotent. It follows from Theorem 2.1 and Corollary 2.4 in \cite{10} that there exists a locally nilpotent derivation $D$ such that $\operatorname{Im}\delta=\operatorname{Im}D$, where $D=(x_2-\frac{1}{2}x_3)\partial_1+x_3\partial_2$. It follows from Theorem 3.4 in \cite{8} that $\operatorname{Im}D$ is a Mathieu-Zhao space of $K[x_1,x_2,x_3]$. Thus, $\operatorname{Im}\delta$ is a Mathieu-Zhao space of $K[x_1,x_2,x_3]$.
\end{rem}

{\bf{Acknowledgement}}:  The second author is very grateful to professor Wenhua Zhao for personal communications about the Mathieu-Zhao spaces. She is also grateful to the Department of Mathematics of Illinois State University, where this paper was partially finished, for hospitality during her stay as a visiting scholar.

\end{document}